\documentclass[EJP]{ejr}

\def\Ib{{\mathbf I}\,}

\def\eopp{{ \vrule height7pt width7pt depth0pt} }
\def\Prf{\noindent {\em Proof:}\, }

\def\Ib{{\mathbf{I}}\,}

\def\xb{{\bf{x}}\,}
\def\thb{{\mathbf{\theta}}\,}
\def\th{{\mathbf{\theta}}\,}

\def\Atl{\tilde{A}\,}
\def\Btl{\tilde{B}\,}
\def\Ctl{\tilde{C}\,}
\def\Ch{\hat{C}\,}

\def\Otl{\tilde{O}\,}

\def\Ctl{\tilde{C}\,}

\def\Vtl{\tilde{V}\,}

\def\Qh{\hat{Q}\,}

\def\eopp{{ \vrule height7pt width7pt depth0pt} }
\def\Prf{\noindent {\em Proof:}\, }

\def\BibTeX{{\rm B\kern-.05em{\sc i\kern-.025em b}\kern-.08em
    T\kern-.1667em\lower.7ex\hbox{E}\kern-.125emX}}

\makeatletter
\@addtoreset{equation}{section}
\makeatother
\newtheorem{thm}{Theorem}[section]
\newtheorem{defn}[thm]{Definition} 
\newtheorem{lem}[thm]{Lemma} 
\newtheorem{cor}[thm]{Corollary} 
\newcommand{\BEQ}{\begin{equation}}
\newcommand{\EEQ}{\end{equation}}
\newcommand{\NEQ}{\end{equation}}

\SHORTTITLE{Hessenberg Input Normal Representations}
\TITLE{Hessenberg Input Normal Representations}
\AUTHORS{Kurt S.~Riedel\footnote{Millennium Partners LLC, New York
    \EMAIL{kurt.riedel@gmail.com}}}

\VOLUME{0}
\YEAR{2005}
\PAPERNUM{0}
\DOI{EDICS 2-SYSM, 2-IDENT}

\ABSTRACT{
Every stable controllable input pair $(\Atl,\Btl)$ is equivalent to
an input pair which is in Hessenberg form and is input normal
(${ A} { A}^* + { B} {B}^* =  { \Ib} $). $(A,B)$ is represented as a submatrix
of the minimal number of Givens rotations. The representation is shown to be generically identifiable.
This canonical form allows for fast state vector
updates and improved conditioning of system identification problems.
}

\AMSSUBJ{93B10, 93B11 93B30}

\begin{document}
\date{February 5, 2005}
\section{INPUT NORMAL PAIRS} \label{I}

The goal of this paper is to propose representations of system pairs $(A,B)$ that are both
well conditioned for system identification and are fast and convenient for
numerical computation. Here the advance matrix, $A$, is a real $n\times n$ matrix
and the control matrix, $B$, is a real $n\times d$ matrix.
The computational advantages of Hessenberg form for state space systems are well known.
\cite{LaubL,VDV}. Recently, studies have shown that input normal form results in significant improvement in 
the conditioning of system identification and filter design problems \cite{MRbook,MR1,MR4}.
In this paper, we show that every stable controllable input pair are equivalent to an input pair with
both properties. By placing the additional requirement that $A_{i+1,i} \ge 0$, this Hessenberg
input normal (HIN) form is generically unique.

We  parameterize the representation using Givens rotations.
Only $nd$ Givens rotations are needed, so state vector advances take ${\cal O}(4nd)$
multiplications. Furthermore the mapping between the Givens angles and the HIN forms is generically
unique. Thus this representation results in generically identifiable problems of system identification.

Let $(B|A)$ denote the concatenation of $B$ with $A$.
Our parameterizations of the $(A,B)$ by 
the  product of Givens matrices is of the form:
\BEQ \label{GenForm}
\left(  B  \ | \ A  \right) =
\left[ \prod_{i=1}^{nd}
G_{j(i),k(i)}(\theta_{i})\right]_{1:n,1:(n+d)} \ \ .
\NEQ
Here $G_{j,k}$ is a Given's rotation in the $i$th and $j$th coordinate by $G_{ij}$:
$g_{i,i}=g_{j,j}=\cos(\theta)$, $g_{i,j}= -g_{j,i}= \sin(\theta)$ and
$g_{k,m} = \delta_{k,m}$ otherwise, where $g_{k,m}$ are the elements of
$G_{ij}$.

Let $(A,B)$ be stable and controllable. We define the controllability Grammian,
$P_{A,B}$ by
\BEQ \label{SteinEq}
P_{A,B} -AP_{A,B}A^*=BB^* \ .
\NEQ
Two systems pairs, $(A, B)$ and  $(\Atl, \Btl)$ are equivalent when
$\tilde{ A} \equiv { T^{-1}AT}$, 
and  $\tilde{ B}\equiv { T^{-1}B}$
for some invertible $T$. If $T$ is orthonormal, we say they are orthogonally equivalent.
We use the freedom of choosing $T$ to choose an equivalent representation of $(A,B)$ which satisfies
\BEQ \label{INeq}
{ A} { A}^* \ + \  { B} {B}^*=  { \Ib_n}  \ \ \ ,
\EEQ
where $\Ib_n$ is $n \times n$ identity matrix.
If $(A,B)$ satisfies (\ref{INeq}), we say $(A,B)$ is {\em input normal} \cite{MR1,MR2}.
If $(A,B)$ is stable and controllable, (\ref{INeq}) is equivalent to requiring
that the controllability Grammian, $P_{A,B}$, equal the identity matrix.

Input normal pairs have superior conditioning and roundoff properties. 
As an  consider the common problem of estimating the observability matrix, $C$
given a stable, controllable input pair, $(A,B)$ and a set of noisy measurements, $\{y_t\}$.
The state vector evolves according to $z_{t+1} = A z_t + B \epsilon_t$,
where $\epsilon_t$ are known random variables. 
The standard recursive least squares (RLS) estimate of $C$ calculates
 $\hat{P}_t$ and $\hat{d}_t$:
\begin{equation}
\hat{P}_t\equiv \frac{1}{t}\sum_{i=1}^t {z}_i{z}_i^{*} , \ \ \
\hat{d}_t= \frac{1}{t}\sum_{i=i}^t  {z}_i {y}_i^* \ .
\end{equation}
The unknown coefficients, $\hat{C}\,$, are estimated by solving
$\hat{P}_t {\hat{C}}=\hat{d}_t$ \cite{LS,MR1}. 
As time increases, $\hat{P}_t$ converges to the controllability Grammian \cite{AM,MR1}
when  the forcing noise, $\epsilon_t$, is white.
When $(A,B)$ is input normal, 
$\hat{P}_t \stackrel{t\rightarrow\infty}{\longrightarrow}{\rm constant} \times \Ib_n$.
Thus the regression estimate of $C$ is well-conditioned.
For more complete analysis
of conditioning in system identification, we refer the reader to \cite{MR4}.

Similarly, IN filter structures are resistant to roundoff error \cite{MRbook}.
Note that input normal  representations will time asymptotically satisfy
the ansatz need by least mean squares (LMS) identification algorithms .
This leads to a second advantage of IN filters: Gradient algorithms such
as the least mean squares (LMS) algorithm often perform well enough in certain
applications to obviate the need for more complicated
and computationally intensive RLS algorithms.

Our representations include the banded orthogonal filters of \cite{MRbook}
as a special case of $d=1$. 
Together with A.~ Mullhaupt \cite{MR5}, we have proposed several other representations of input normal pairs,
where $(A,B)$ is the product of $nd$ Givens rotations.

Our approach has similarities to that of embedded lossless systems (ELS) \cite{Desai,Vaid,VeenVib}.
While embedded lossless systems are of interest in a few applications, we believe that it is simpler to parameterize
$(A,B)$ as a HIN pair. This allows $(C,D)$ to be determined by linear or pseudolinear regression, 
and this regression is very well conditioned 
since the expected value of  controllability Grammian is the identity matrix.
In contrast, ELS satistfy ${ A} { A}^* \ + \  { B} {B}^*+  B_{ext} {B_{ext}}^*=  { \Ib_n}$ where
$B_{ext}$ is an artificial term for the embedding. Thus the ELS controllability Grammian for the actual system pair can be ill-conditioned. 
The construction of the Givens matrix representation of \cite{VeenVib} parallels our representation in Section \ref{HINRepSect}.
One can interpret our results as a simplification of \cite{VeenVib} for those who prefer input normal representations to embedded lossless systems.
Another difference between our treatment and the analyses of \cite{Desai,Vaid,VeenVib} is that
we describe when redundant representations can occur.

{\em Notation:}
By $A_{i:j,k:m}$, we denote the $(j-i+1)\times (m-k+1)$ subblock
of $A$ from row $i$ to row $j$ and from column $k$ to column $m$.
The direct sum of matrices is denoted by $\oplus$.
We denote the matrix transpose of $A$ by $A^*$ with no complex conjugation
since we are interested in the real system case.

\section{Hessenberg Input Normal Pairs}
\label{DefSect}

Hessenberg form is a canonical form
where $A$ is Hessenberg with  the additional restriction
that $B_{1,1} \ge 0$, $B_{j,1}=0$ for $j>1$.

\begin{defn} \label{DefHOess}
The input pair (A,B) is in  Hessenberg form
if $A$ is a Hessenberg matrix, $B_{1,1} \ge 0$, 
and $B_{j,1} =0$ for $j>1$.
A Hessenber pair is nondegenerate if $|B_{1,1}|<1$.
A Hessenberg pair is unreduced
if $A_{i+1,i} \ne 0 $ for $1 \le i <n$ and $B_{1,1} \ne 0$.
A Hessenberg pair is standard if 
$A_{i+1,i} \ge 0 $ for $1 \le i <n$, $0 \le B_{1,1} \le 1 $. 
A Hessenberg pair is strict if it is unreduced and standard.
\end{defn}

If $(A,B)$ is a HIN pair, so is $(EAE^{-1},EB)$ where $E$ is an arbitrary
siganture matrix: $|E_{i,j}| =\delta_{i-j}$.
Thus we seek a representation of standard HIN pairs to eliminate the multiplicity of equivalent 
representations when $(A,B)$ is unreduced.



\begin{thm}\label{HOCon} \cite{LaubL,VDV}
Any observable input pair is {\rm orthogonally} equivalent to
a system in Hessenberg form.
\end{thm}

The standard proof of Theorem \ref{HOCon} begins by transforming
$B$ to its desired form and then defines Givens
rotations which zero out particular elements in $A$ in
successive rows or columns \cite{GVL,LaubL}.

For nondegenerate HIN pairs, we find that the set of standard
input normal pairs has a bijective representation as an easy to parameterize subset
of  Givens product representations.
The existence and uniqueness results in this section originated in the unpublished technical report \cite{MR5}. 
The proof of Theorem \ref{MONthm} is also in \cite{VeenVib}.
We reproduce them for completeness. 

\begin{thm}\cite{MR5}
\label{MONthm}
Every  stable, real controllable output pair
$(A, B)$, is similar to a standard HIN pair.
\end{thm}

\Prf
The unique solution, $P_{A,B}$, of Stein equation, (\ref{SteinEq}),
is strictly positive definite. 
Let $L$ be the  unique Cholesky lower triangular factor of ${ P}$
with positive diagonal entries:
${ P} ={ LL}^{*}$.
We set  $T=L^{-*}$.
Let $U$ be orthogonal transformation that takes $(T^{-1}AT, T^{-1}B)$
to the Hessenberg form as described in \cite{LaubL,VDV}. We now choose the signature matrix, $E$
such that $T\equiv EUL^{-*}$ makes $(T^{-1}AT, T^{-1}B)$ a standard HIN pair.
\eopp


Degenerate HIN pairs correspond to the direct sum of an identity
matrix and a nondegenerate HIN system:

\begin{lem}
Every  stable, real controllable input pair
$(A, B)$, is similar to a HIN pair which is the direct sum 
of the identity matrix and a  nondegenerate HIN pair:
$(B|A) = \left(\Ib_m  \oplus \Qh\right)$
for some $m \le d$, where $\Qh$ is a $(n-m) \times (n+d-m)$row orthogonal matrix.
\end{lem}


\label{UniqOvSect} \label{HINUniqSect}
\label{UniqSect}


There are two main ways in which one of our system representations can fail to
parameterize linear time invariant systems in a bijective fashion.
First, there may be a multiplicity of equivalent HIN systems.
Second, Givens product representation
such as (\ref{GenForm}) may have multiple (or no) parameterizations
of the same input pair.
We now show when that if $(A,B)$ is a strict HIN pair,there are no equivalent strict HIN pairs. 
The uniqueness results are based on the following
lemma that generalizes the Implicit Q theorem \cite{GVL}
to HIN pairs:

\begin{lem} \label{HoQ}
Let $(A, B)$ and $(\Atl,\Btl)$ be orthogonally 
equivalent standard nondegenerate HIN pairs
(${\Atl} \equiv { T^{-1}AT}, { \Btl}\equiv { T^{-1}B}$ where $T$ is orthogonal).
Let $A_{k+1, k}= 0$, $B_{1,1}>0 $ and  $A_{j+1,j}>0$ for $j<k$,
then $T= \Ib_{k} \oplus U_{n-k}$, where $ U_{n-k}$ is an
$(n-k) \times (n-k)$ orthogonal matrix. Furthermore, $k>1$
and $\Atl_{k+1, k}= 0$.
\end{lem}

Since $B_{j,1}= \Btl_{j,1}=0 $, $j>1$, $T_{j,1}= \delta_{j,1}$.
The result follows from  the Implicit Q theorem \cite{GVL}.
\eopp

\begin{cor} \label{HoUniq}
If $(A, B)$ is a strict nondegenerate HIN pair, then there are no
other equivalent strict HIN pairs.
\end{cor}

When $A_{k+1, k}= 0$, there are many different equivalent HIN pairs. This degeneracy may be 
reduced or eliminated if additional restrictions are imposed on $(A,B)$. In \cite{MR5}, it is proposed that
$A_{k+1:n,k+1:n}$ be placed in real Schur form with a given ordering of the eigenvalues and standardization
of the $2 \times 2$ diagonal subblocks. An alternative restriction is to require $B_{k+1,2}\ge 0$ and $B_{k+2:n,2}=0$
when  $d>1$. 

\section{Givens Representation of Hessenberg input normal form.} \label{HINRepSect}

In this section, we give representation results for HIN pairs.
The first column of $B$ satisfies $0 < B_{1,1} \le 1 $,
$B_{j,1} = 0$ for $j>1$. 
We use Givens rotations to zero out  the lower diagonal of $A$ and the row $2$ through
row $d$ of $B$. For each row of the $(B | A)$, we use $d$ Givens rotaions.
to zero out the lower diagonal of $A$ and columns 2 through $d$ of $B$. 
The following lemma shows how this works on a single row of $(B|A)$.

\begin{lem} \label{Giv}
Let ${X}$ be the set of real $d+1$ tuples, $\xb= (x_1,\ldots x_d, x_{d+1})$ 
where $\sum x_j^2 =1$ and $x_d \ge 0$. 
Let $g(\thb)$ denote the $d+1$ tuple 
$(0 \ldots 0, 1)G_{d+1,d}(\theta_d)G_{d,d-1}(\theta_{d-1})\ldots  G_{3,2}(\theta_2) G_{2,1}(\theta_1)$,
where $G_{i,i+1}$ are the Givens rotations in $R^{d+1}$.
Let $\Theta$ denote the restriction of the Givens rotation angles
to $0 \le \theta_1 \le \pi$, $-\pi/2 < \theta_j \le \pi/2$,
for $1<j \le d$. The mapping from the $\Theta$, the domain of the Givens angles, is \underline{onto} $X$.
The mapping is \underline{one-to-one} for the set of $\xb$ where $|x_1|>0$.
\end{lem}

\Prf
Clearly, the first Givens angle is uniquely determined by the last component of $\xb$. If $x_{d+1}= \pm 1$, 
the rest of the representation, $\{\theta_2 \ldots \theta_d\}$ is arbitrary. 
The $k$th component of $g(\thb)$ is $x_k=g(\thb)_k = \cos(\theta_k) \prod_{i=k+1}^d \sin(\theta_i)$.
If $x_{k-1}$ is nonzero, we determine $\theta_k$ uniquely by $g(\thb)_k =x_k$ and 
$sign(g(\thb)_k) = sign(x_{k-1})$. If  $x_{k-1}=0$, we require $sign(g(\thb)_k) = sign(x_{j})$ 
where $j$ is the largest index, $1 \le j < k$
with $x_j \ne 0$. This determination of the Givens angles is unique unless $x_1$ vanishes.
\eopp

Note that $\theta_j = -\pi/2$ is not necessary.
We can make the mapping of $g(\thb): \Theta \rightarrow X$ globally one to one if we require that whenever $\sum_{i<j} |x_i| =0$ for some $j$,
that $\theta_i= 0$ for $i<j$ and that if $j<d$ then $\theta_j=0$. If $j=d$, then $\theta_j = \pm \pi$.

\begin{thm} \label{HINRepthm}
Every real  HIN pair has the representation:
\BEQ \label{HINRepEq}
(B|A) = (0_{n,d} | \Ib_n )U^{(1)}(\theta_1:\theta_d) V^{(2)}(\theta_{d+1}:\theta_{2d})  \ldots
V^{(n-1)} V^{(n)}(\theta_{(n-1)d+1}:\theta_{nd})
\NEQ
where $0_{n,d}$ is the $n \times d$ matrix of zeros.
Here $U^{(1)}$ and $V^{(k)}$ are $(n+d) \times (n+d) $ matrices which are the product of $d$ Givens rotations:
\BEQ \label{VkStruct}
V^{(k)}(\theta_{(k-1)d+1}:\theta_{kd} ) = \\
G_{k+d,k+d-1}(\theta_{(k-1)d+1}) G_{n+d-1,d}(\theta_{(k-1)d+2}) \ldots G_{4,3}(\theta_{kd-1}) G_{3,2}(\theta_{kd})
\NEQ
\BEQ \label{U1Struct}
U^{(1)}(\th_1:\th_d) = \\
G_{d+1,1}(\theta_1) G_{1,2}(\theta_2)  \ldots  G_{d-2,d-1}(\theta_{d-1}) G_{d-1,d}(\theta_d)  
\NEQ
Every matrix of the form of the righthand side of (\ref{HINRepEq}) is a HIN pair. 
\end{thm}

We  successively determine the Givens angles, $\theta_k$.
At the $k$th stage, $\theta_{(n-k-1)d+1} \ldots \theta_{(n-k)d}$ are determined
to zero out the $d$ of the $d+1$ nonzero entries
in the $k$th row.
By orthogonality, the other entries in the $n-k$th column must be zero.

\Prf
Let
\BEQ
\Gamma^{(k)}(\theta_{dk+1}: \theta_{nd} ) \equiv (B|A)
V^{(n)*}  V^{(n-1)*} \ldots  V^{(k+1)*} \ .
\NEQ
Assume that $\Gamma^{(k)}$
has its last $k$ rows  satisfying
$\Gamma^{(k)}_{i,j} = \delta_{i-j+d}$.
Since  $\Gamma^{(k)}$ has orthonormal rows, the last $k$ columns satisfy
$\Gamma^{(k)}_{i,j}= \delta_{i-j+n}$.
Select $\theta_{kd+1}$ through $\theta_{(k+1)d}$ such
$\Gamma^{(k+1)}_{1:d,k+1} =0$ and $\Gamma^{(k+1)}_{n-k+d-2,n-k-1} =0$. Then
$\Gamma^{(k+1)}_{i,j} =\delta_{i-j+d}$ for the last $k+1$ rows and the last $k+1$ columns.

To show all Givens products of the form in (\ref{HINRepEq}) generate HIN pairs, consider
\BEQ \label{Xdef}
X^{(k)} (\theta_1 \ldots  \theta_{kd})=(0 | \Ib_n )U^{(1)}(\theta_1:\theta_d) V^{(2)}(\theta_{d+1}:\theta_{2d})  \ldots
V^{(k)}(\theta_{(k-1)d+1}: \theta_{kd}) 
\NEQ
$X^{(k)}$ is in Hessenberg form for each $k$.
\eopp

For $d=1$ and $B_{1,:}=0$, (\ref{HINRepEq}) is the well-known expression of an
unitary Hessenberg
matrix as a product of $n$ Givens rotations \cite{AGR}.
The fast filtering methods of \cite{VDV,MR1} may be further sped up when
 $(B|A)$  is a submatrix product of $nd$ Givens rotations.

\begin{thm}\label{One}
Let $\Theta$ be the set $0 \le \theta_{kd+1} \le \pi$, $-\pi/2 < \theta_{kd+j} \le \pi/2$, 
where $1<j \le d$. 
The representationsof  Theorem \ref{HINRepthm} maps $\Theta$ \underline{onto} the set of standard HIN pairs. 
The mapping is \underline{one-to-one} for $d=1$. For $d>1$, the mapping is \underline{one-to-one} on the set where $B_{1,d}>0$, $B_{j,2}>0$ for $j>1$.
Let $B^{cat}$ be the concatenation of $B_{:,2:d}$ with the vector $(B_{1,1,},A_{2,1}\ldots A_{n,n-1})$.
If whenever $\sum_{i<j} |B^{cat}_{k,i}| =0$, we impose the constraint that $\theta_{(k-1)d+i}= 0$ for $i<j$ and for $j<d$ the additional 
constraint that $\theta_{(k-1)d+j}>0$, 
then there is a one to one correspondence between strict HIN pairs and
the parameterization of Theorem \ref{HINRepthm}
\end{thm}

\Prf
We repeatedly apply Lemma \ref{Giv} to each row of $X^{(k)}$.
\eopp

Theorem \ref{One} does not address the multiplicity of equivalent HIN pairs when $A_{i+1,i}$ vanishes.
The representation of Theorem \ref{HINRepthm} is nonuniform because the order of the Givens rotations in $U^{(1)}$ differ from that of $V^{(k)}$.
If we require $B_{1,d}\ge 0$ and $B_{j,d}=0$ instead of $B_{1,1}\ge 0$ and $B_{j,1}=0$, the representation is more uniform
in selecting the order of the Givens rotations.

The disadvantage of the Givens angle restrictions is that the representation is discontinuous when $(A,B)$ corresponds to a boundary value 
of $\Theta$. An example of this behavior is that if $B_{n,d}$ is near zero, $\theta_1$ varies from $\pi/2$ to $-\pi/2$.
For this reason, it may be preferable to eliminate the angular restrictions on $\thb$ 
while performing a numerical optimization to determine the Givens angles.

\section{Summary.}

We have examined the uniqueness/identifiability of system pairs when $(A,B)$ is simultaneously in Hessenberg form and input normal.
Since controllability Grammian is the identity matrix, estimation of the observability matrix is well-conditioned.
These system pairs have good roundoff properties.
To transform a specific input pair to IN form, the Stein equation (\ref{SteinEq}) must be solved
for $P_{A,B}$.
The numerical conditioning of this problem can be quite poor
\cite{PenzlBnd,MR4}. Thus it is much better to start with an explicit representation of $(A,B)$ in HIN form 
than to transform a non-input normal system into HIN form.
To make IN pair attractive for applications, we describe a concise representation with fast state vector updates.

In this case, the representation of Theorem \ref{HINRepthm} is particularly convenient.
$(A,B)$  has a representation as a submatrix of the product of $nd$ Givens matrices.
We have shown how to place restrictions on the parameters such that
the Givens representaion is in one to one correspondence
with the set of standard HIN pairs.
For statistical estimation,
it is sometimes desirable to eliminate redundancy
in the parameterization.
We address redundancy in two ways. First, we categorize
when two distinct HIN pairs are equivalent.
Second, we impose constraints on
the parameters in (\ref{GenForm}) to eliminate redundant parameterizations
of the same HIN pair \underline{generically}.
However the resulting constraint set, $\Theta$ has the property that small perturbations of $(A,B)$ can cause large
changes in the Givens angles near the boundary of $\Theta$. Thus it may be better to not restrict the Givens angles in
numerical optimizations.

We do not explicitly store or multiply by $V^{(k)}$ in (\ref{HINRepEq}).
Instead we store only the Givens parameters and
we perform the matrix multiplication implicitly. 
Thus vector multiplication by $A$ and by $ \frac{d}{d \theta_k} A$ require $O(4nd)$ and $O(8nd)$ operations.
Similar Givens product representations of output normal $(A,C)$ when $A$
is in real Schur form or (A,C) or in Hessenberg observer form may be found in  
is in \cite{MR5}.

In summary, HIN representations offer the best possible
conditioning while having a convenient representation with fast
matrix multiplication.

\bibliographystyle{IEEE}

\ACKNO{The author thanks Andrew Mullhaupt for a long and fruitful collaboration in which many related
results were discovered.}
\end{document}